\documentclass[13pt]{amsart}
\usepackage{amsmath,amsfonts,amssymb,amscd}
\usepackage[all]{xy}
\usepackage{graphics}
\usepackage{epsfig}
\usepackage{psfrag}

\def\ignore #1 {}

\newtheorem{ass}{Assumption}
\newtheorem{lem}{Lemma}
\newtheorem{thm}{Theorem}

\newtheorem{dfn}{Definition}
\newtheorem{cor}{Corollary}
\newtheorem{rmk}{Remark}

\newtheorem{ex}{Example}

\begin{document}

\title{Equivariant Gerbes on Complex Tori}
\author{Oren Ben-Bassat}
 \maketitle

\centerline{\bf }

%%%%%%%%%%%%%%%%%%%%%%%%%%%%%%%%%%%%%%%%%%%%%
%%%%%%%%%%%%%%%%%%%%%%%%%%%%%%%%%%%%%%%%%%%%%
%%%%%%%%%%%%%%%%%%%%%%%%%%%%%%%%%%%%%%%%%%%%%

%%%%%%%%%%%%%%%%%%%%%%%%%%%%%%%%%%%%%%%%%%%%%
%%%%%%%%%%%%%%%%%%%%%%%%%%%%%%%%%%%%%%%%%%%%%
%%%%%%%%%%%%%%%%%%%%%%%%%%%%%%%%%%%%%%%%%%%%%

%%%%%%%%%%%%%%%%%%%%%%%%%%%%%%%%%%%%%%%%%%%%%
%%%%%%%%%%%%%%%%%%%%%%%%%%%%%%%%%%%%%%%%%%%%%
%%%%%%%%%%%%%%%%%%%%%%%%%%%%%%%%%%%%%%%%%%%%%

%%%%%%%%%%%%%%%%%%%%%%%%%%%%%%%%%%%%%%%%%%%%%
%%%%%%%%%%%%%%%%%%%%%%%%%%%%%%%%%%%%%%%%%%%%%
%%%%%%%%%%%%%%%%%%%%%%%%%%%%%%%%%%%%%%%%%%%%%

\begin{abstract}  We explore a new direction in representation theory which comes from holomorphic gerbes on complex tori.  The analogue of the theta group of a holomorphic line bundle on a (compact) complex torus is developed for gerbes in place of line bundles.  The theta group of symmetries of the gerbe has the structure of a Picard groupoid.  We calculate it explicitly as a central extension of the group of symmetries of the gerbe by the Picard groupoid of the underlying complex torus.  We discuss obstruction to equivariance and give an example of a group of symmetries of a gerbe with respect to which the gerbe cannot be equivariant.   We calculate the obstructions to invariant gerbes for some group of translations of a torus to be equivariant.    We survey various types of representations of the group of symmetries of a gerbe on the stack of sheaves of modules  on the gerbe and the associated abelian category of sheaves on the gerbe (twisted sheaves).  
\end{abstract}

\section{Introduction}\label{intro}  

This article is a continuation of \cite{BGerbes} so let us briefly summarize some of the pertinent points of that work.  The symbols $\mathcal{O}$ and $\mathcal{O}^{\times}$ refer to the holomorphic and nowhere zero holomorphic functions on some complex manifold (the specific manifold will be clear from the context) with the obvious group structures.  Holomorphic (banded) $\mathcal{O}^{\times}$-gerbes on a complex manifold $M$ form a $2-$groupoid and their equivalence classes are calculated by the cohomology $H^{2}(M,\mathcal{O}^{\times})$ computed in the classical topology.

A $\mathcal{O}^{\times}${\it -gerbe} $\mathfrak{G}$ over a complex manifold $M$ is a twisted form of $\mathcal{P}ic_{M}$.  Let us give a quick description of an $\mathcal{O}^{\times}$-gerbe without all the details.  It consists first of all of a stack $\mathfrak{G}$ on $M$.  This includes a category $\mathfrak{G}(U)$ for each open set $U \subset M$, and functors $\mathfrak{G}(U) \to \mathfrak{G}(V)$ for every inclusion of open sets $V \subset U$ satisfying certain axioms.  Second, one has an action functor of the stack in Picard groupoids $\mathcal{P}ic_{M}$ on $\mathfrak{G}$.
\[\mathcal{P}ic_{M} \times \mathfrak{G} \to \mathfrak{G}
\]
\[(\mathcal{L},G) \mapsto \mathcal{L} \cdot G
\]
such that for every $m \in M$ there is an open $U$ containing $m$ such that $\mathfrak{G}(U)$ is non-empty and for each object $G$ of $\mathfrak{G}(U)$, the functor 
\[\mathcal{P}ic_{M}(U) \to \mathfrak{G}(U)
\]
\[\mathcal{L} \mapsto \mathcal{L} \cdot G
\]
is an equivalence of categories.  $\mathcal{O}^{\times}$-gerbes on $M$ can be tensored over $\mathcal{P}ic_{M}$ giving the structure of a group to the set of equivalence classes which agrees with the cup product on cohomology classes.
Some material on holomorphic gerbes both old and recent can be found in  \cite{Chatterjee}, \cite{Gi1971}, \cite{Br1993}, \cite{HL}, \cite{Hao}, \cite{FeHeRoZh2008}, \cite{BGerbes}, \cite{BB}, \cite{Po2008}, \cite{Y2010} and \cite{Hi2010}.   In this paper, we will refer to  $\mathcal{O}^{\times}$-gerbes simply as grebes.  Gerbes arise in physics in attempting (see for instance \cite{GKGG}) to understand ambiguities in the definition of the generalized K\"{a}hler potential in the context of sigma models with $N=(2,2)$ supersymmetry.  

Let $X=V/\Lambda$ be a complex torus.  The projection map will be denoted by 
\[\pi:V \to X. 
\]
Given an element $E \in \text{Alt}^{3}(\Lambda,\mathbb{Z})$ such that  
\begin{equation}\label{str1}
E(x,y,z) = E(ix,iy,z)+E(x,iy,iz)+E(ix,y,iz)
\end{equation} 
we defined an element of $C^{2}(\Lambda,\mathcal{O}(V))$ by 
\begin{equation}\label{eqn:Heqn}
\begin{split}
H_{\lambda_{1},\lambda_{2}}(v) &= \frac{1}{8}\Big{(}E(v,\lambda_{1},\lambda_{2})+\frac{1}{2}E(iv,i\lambda_{1},\lambda_{2})
+\frac{1}{2}E(iv,\lambda_{1},i\lambda_{2})\Big{)} \\
&  \quad +\frac{i}{8}\Big{(}\frac{1}{2}E(v,i\lambda_{1},\lambda_{2})+\frac{1}{2}E(v,\lambda_{1},i\lambda_{2})- E(iv,\lambda_{1},\lambda_{2})\Big{)}.
\end{split}
\end{equation}  If we are also given
an element $B \in \text{Alt}^{2}(\Lambda, \mathbb{R})$ we defined 
 $\Phi= \Phi^{(B,E)}$ of $C^{2}(\Lambda, \mathcal{O}^{\times}(V))$ by
\begin{equation}\label{TheCocycle}\Phi^{(B,E)}_{\lambda_{1}, \lambda_{2}}(v)=\exp\Big{(} \frac{1}{2}B(\lambda_{1},\lambda_{2})+ H_{\lambda_{1},\lambda_{2}}(v) + \beta'_{\lambda_{1},\lambda_{2}}+i\beta''_{\lambda_{1},\lambda_{2}}\Big{)}
\end{equation} 
which in fact is shown to lie in $Z^{2}(\Lambda, \mathcal{O}^{\times}(V))$.  The definition of $\beta', \beta'' \in C^{2}(\Lambda,\mathbb{R})$ is not important for us and can be found in \cite{BGerbes}.  By applying to $\Phi$ the connecting homomorphism in the short exact sequence 
\[0 \to \mathbb{Z} \to \mathcal{O} \to \mathcal{O}^{\times} \to 1
\]
we get a class in $Z^3(\Lambda,\mathbb{Z})$ the skew-symmetrization of which is an element of $\text{Alt}^{3}(\Lambda, \mathbb{Z})$  which agrees with $E$.  If $E=0$ then the element $\Phi^{B,E}$ comes from the image of $B$ in $Z^{2}(\Lambda,\mathcal{O}(V))$.  The map $\text{exp}$ is defined by $\text{exp}(z)= e^{2\pi i z}$.  As shown in \cite{BGerbes}, all classes in 
\[H^{2}(X,\mathcal{O}^{\times}) \cong H^{2}(\Lambda, \mathcal{O}^{\times}(V))\] come from such a pair $(B,E)$.  The cohomology on the right is the group cohomology of the lattice $\Lambda$ acting on functions by translation.  The class $E$ is uniquely determined although the class $B$ is not.  In fact two pairs $(B_{1},E_{1})$ and $(B_{2}, E_{2})$ describe isomorphic gerbes if and only if $E_{1}=E_{2}$ and $B_{1}- B_{2} \in \text{Alt}^{2}(\Lambda, \mathbb{Z})+ \text{Alt}^{2}(\Lambda, \mathbb{R})^{(1,1)}$ where the space $\text{Alt}^{2}(\Lambda, \mathbb{R})^{(1,1)}$ depends on the complex structure on the complex torus $X$. This is a ``canonical" cocycle in the same sense as the Appell-Humbert theorem \cite{Ap1891},\cite{Hu1893}, \cite{We1958}.  
A {\it $K-$equivariant gerbe} on a complex manifold $M$ acted on by some group $K$ consists of a gerbe $\mathfrak{G}$ on $M$, a $(1-)$isomorphism of the two pullbacks of $\mathfrak{G}$ to $K \times M$ (which is analytic over each point of $K$), and a $(2-)$ isomorphism relating the three resulting (1-)isomorphisms on $ K \times K \times M$ (analytic over each pair of points in $K$) which satisfies a natural coherence condition on $K \times K \times K \times M$.  A {\it $K$-invariant gerbe} is simply a gerbe $\mathfrak{G}$ on $M$, a $(1-)$isomorphism of the two pullbacks of $\mathfrak{G}$ to $K \times M$  which is analytic over each point of $K$. 

In this paper we address the question of when a $K$-invariant gerbe has a $K$-equivariant structure.  $K$-invariant gerbes are classified by $H^{2}(M,\mathcal{O}^{\times})^{K}$.  Geometrically, we compute obstructions to a $K$-invariant gerbe on $M$ being the pullback of a gerbe on a quotient stack $[M/K],$ which is the same as an equivariant gerbe.  We conclude the introduction with two technical remarks that will be used in the body of this article.

\begin{rmk}
There is a spectral sequence converging to the equivalence classes of equivariant gerbes where we take in the following group cohomology with differentiable cochains 
\[K^{i} \to H^{j}(M,\mathcal{O}^{\times})
\]
(for values of $i$ and $j$ shown below).  We can also replace $K$ by one of its discrete subgroups.

The $E_{2}$ term of this spectral sequence looks in the lower left corner like
\def\csg#1{\save[].[dddrrr]!C*+<6.5pc,0pc>[F-,]\frm{}\restore}
\begin{equation}\label{ss}
\xymatrix@=.7pc{
 & \csg1
&
 &
 &   
 \\
 &  H^{2}(M,\mathcal{O}^{\times})^{K} &
 &
&    \\
 & H^{1}(M,\mathcal{O}^{\times})^{K}
&   H^{1}(K,H^{1}(M,\mathcal{O}^{\times}))&
H^{2}(K,H^{1}(M,\mathcal{O}^{\times})) &   \\
 & (\mathcal{O}^{\times}(M))^{K} & H^{1}(K, \mathcal{O}^{\times}(M)) &H^{2}(K, \mathcal{O}^{\times}(M))  & H^{3}(K, \mathcal{O}^{\times}(M))   \\
& &  &  & 
}.
\end{equation}
converging along the third anti-diagonal to the equivalence classes of equivariant gerbes.  Thus there are two maps that concern us:
\[H^{2}(M,\mathcal{O}^{\times})^{K} \stackrel{d_{2}}\to H^{2}(K,H^{1}(M,\mathcal{O}^{\times}))
\]
and 
\[H^{2}(M,\mathcal{O}^{\times})^{K} \supset \text{ker}(d_{2}) \stackrel{d_{3}}\to H^{3}(K, \mathcal{O}^{\times}(M))/d_{2}(H^{1}(K,H^{1}(M,\mathcal{O}^{\times})))
.\]
In order to answer the question of  when a $K$-invariant gerbe has a $K$-equivariant structure we will therefore need to characterize the conditions under which an element of the $E_{2}$ term concentrated in $H^{2}(M,\mathcal{O}^{\times})^{K}$ survives to the third anti-diagonal of $E_{4}$ which is the third anti-diagonal of  $E_{\infty}$.  Notice also that the above spectral sequence is contravariant with respect to maps of pairs $(K_{1},M_{1}) \to (K_{2},M_{2})$ which respect the group and analytic space structures and are compatible with the group action. 
\end{rmk}
\begin{rmk}\label{annoying}
Suppose now that $X=M$ is a complex torus and $K$ is a torsion subgroup of the translations of $X$ such that every $x \in K$ fixes a gerbe $\mathfrak{G}$.  If we consider the short exact sequence 
\[0 \to \text{Pic}^{0}(X) \to \text{Pic}(X) \to \text{NS}(X) \to 0
\]
then there is a resulting long exact sequence: 
\[\cdots \to [H^{1}(K, \text{NS}(X)) = \text{Hom}(K,\text{NS}(X))] \to H^{2}(K,   \text{Pic}^{0}(X)) \to H^{2}(K,   \text{Pic}(X)).
\]
Since $K$ is torsion the term $\text{Hom}(K, \text{NS}(X))$ vanishes and so $H^{2}(K,   \text{Pic}^{0}(X))$ includes into  $H^{2}(K,   \text{Pic}(X))$.  Also, we can conclude that any element in $H^{1}(K, \text{Pic}(X))$ can be lifted to $H^{1}(K, \text{Pic}^{0}(X))= Hom(K, \text{Pic}^{0}(X))$.  Since any such homomorphism can clearly be lifted to a homomorphism from $K$ to $Z^{1}(\Lambda, U(1))$, the differential $d_{2}$ vanishes on $H^{1}(K,H^{1}(M,\mathcal{O}^{\times}))$.
\end{rmk}
%\begin{ack}
%I would like to thanks R. Donagi, T. Pantev, and C. Daenzer for helpful comments.
%\end{ack}
\section{Generalities}\label{Generalities}
\begin{dfn}\label{Kdef}Let $\mathfrak{G}$ be a gerbe on complex torus $X$.  We define $K(\mathfrak{G}) \subset X$ to be the subgroup of $X$ consisting of elements $x \in X$ such that the translation $x \cdot \mathfrak{G}$ of $\mathfrak{G}$ is isomorphic to $\mathfrak{G}$.  The inverse image under $\pi:V \to V/\Lambda$ is denoted by 
\[\Lambda(\mathfrak{G}) = \pi^{-1}(K(\mathfrak{G}))
\]
\end{dfn}
We will often consider images of certain classes under the pullback map on group cohomology
\[\pi^{*}:H^{j}(K(\mathfrak{G}),M) \to H^{j}(\Lambda(\mathfrak{G}),M)
\]
defined for a $K(\mathfrak{G})$-module $M$.   
\begin{dfn}Let $\mathfrak{G}$ be a gerbe on a complex torus $X$.
The theta groupoid $\mathcal{G}(\mathfrak{G})$ is defined as follows
\[\text{ob}(\mathcal{G}(\mathfrak{G})) = \{(x,f)| x \in X, f: x \cdot \mathfrak{G}  \xrightarrow{\cong} \mathfrak{G}\}
\]
For $x_{1} \neq x_{2}$
\[\text{Hom}_{\mathcal{G}(\mathfrak{G})}((x_{1},f_{1}), (x_{2}, f_{2})) = \emptyset
\]
and
\[\text{Hom}_{\mathcal{G}(\mathfrak{G})}((x,f_{1}), (x, f_{2})) = Isom(f_{1}, f_{2}).
\]
If one thinks of $f_{1}$ and $f_{2}$ as twisted line bundles for $(x\cdot \mathfrak{G})^{-1} \otimes \mathfrak{G}$ then $f_{1}^{-1}f_{2}$ is a line bundle on $X$ and an isomorphism is a global nowhere vanishing section of this line bundle.
The Picard structure on this groupoid is given by the functor
\[m: \mathcal{G}(\mathfrak{G}) \times \mathcal{G}(\mathfrak{G}) \to \mathcal{G}(\mathfrak{G})
\]
defined on objects by
\[m((x_{1},f_{1}), (x_{2},f_{2})) = (x_{1}+x_{2}, f_{1} \circ (x_{1} \cdot f_{2}))
\]
and on morphisms in the obvious way.  
Similarly, we have the lifted theta groupoid $G(\mathfrak{G})$.  It is defined in the same way as the theta groupoid but using elements of the vector space $V$ acting via $\pi$ in place of elements of $X$.
\end{dfn} 
Notice that we have a short exact sequences of Picard groupoids
\[0 \to \mathcal{P}ic(X) \to \mathcal{G}(\mathfrak{G}) \to K(\mathfrak{G}) \to 0.
\]
and similarly
\[0 \to \mathcal{P}ic(X) \to G(\mathfrak{G}) \to \Lambda(\mathfrak{G}) \to 0.
\]
corresponding to the action of $K(\mathfrak{G})$ on  $\mathcal{P}ic(X)$ by translations, for which $\mathcal{P}ic(X)$ becomes a $K(\mathfrak{G})$-module.  

If we pass to equivalence classes we get short exact sequences 
\begin{equation}
\label{fristextobs}
0 \to \text{Pic}(X) \to \pi_{0}(\mathcal{G}(\mathfrak{G})) \to K(\mathfrak{G}) \to 0
\end{equation}
and similarly
\begin{equation}
\label{LatticeSES}
0 \to \text{Pic}(X) \to \pi_{0}(G(\mathfrak{G})) \to \Lambda(\mathfrak{G}) \to 0.
\end{equation}
We call $\pi_{0}(\mathcal{G}(\mathfrak{G}))$ and $\pi_{0}(G(\mathfrak{G}))$ the theta group and lifted theta group of $\mathfrak{G}$.
In general the arrows $w \cdot \mathfrak{G} \to \mathfrak{G}$ are twisted line bundles represented by cochains $\tau^{w} \in C^{1}(\Lambda,\mathcal{O}^{\times}(V))$ whose precise formula will be derived in Section \ref{ExplicitChoice}.

\begin{ex}
In the case where $\mathfrak{G}$ is trivial, then $K(\mathfrak{G}) = X$ and the corresponding structure of a Picard groupoid on the product
\[X \times \mathcal{P}ic(X)
\]
is
\begin{equation}\label{SplitGrpStr}
(x_{1},L_{1}) \cdot (x_{2}, L_{2}) = (x_{1}+x_{2}, L_{1} \otimes (x_{1} \cdot L_{2})).
\end{equation}
\end{ex}
The automorphisms of a gerbe on any space form a Picard groupoid which is canonically
isomorphic to the Picard goupoid of $\mathcal{O}^{\times}-$torsors on that space.   Given a gerbe $\mathfrak{G}$ and
elements $x_{1},x_{2} \in K(\mathfrak{G})$, we actually get an $\mathcal{O}^{\times}-$torsor $\mathcal{L}^{x_{1},x_{2}}$ on $X$ defined as the obstruction to lifting translations to the
gerbe in an additive manner.  In other words, we have the following definition
\begin{dfn}\label{LObsDef}
The $\mathcal{O}^{\times}$-torsor  $\mathcal{L}^{x_{1},x_{2}}$ is the automorphism of $\mathfrak{G}$ given by the composition of isomorphisms 
\begin{equation}\label{GerbeSequence}
\mathfrak{G} \to (x_{1} + x_{2}) \cdot \mathfrak{G} = x_{1} \cdot (x_{2} \cdot \mathfrak{G}) \to x_{1} \cdot \mathfrak{G} \to \mathfrak{G}.
\end{equation}
\end{dfn}

We will calculate the pullback to $H^{2}(\Lambda(\mathfrak{G}),Pic(X))$ of the class in
\[H^{2}(K(\mathfrak{G}), \text{Pic}(X)) 
\]
which classifies the extension (\ref{fristextobs}) and corresponds to the image of $[\mathfrak{G}] \in H^{2}(X,\mathcal{O}^{\times})^{K(\mathfrak{G})}$ under the map $d_{2}$ from the spectral sequence (\ref{ss}).

\section{A Criterion For Invariance}\label{Criterion}

In this section we use the canonical cocycle to see how a gerbe $\mathfrak{G}$ defined by the cocycle $\Phi$ defined in  (\ref{TheCocycle})  transforms
under the pullback by a translation.  The main point is that after translating a cocycle, one can put it back in canonical form to see how it changes.  A subgroup of the torus will fix the gerbe
and this subgroup has a canonical central extension.  This extension involves choices of an isomorphism of a gerbe and its translation.

Translation by $w \in V$ acts trivially on $\Lambda=\pi_{1}(X)$ and so it
multiplies the cocycle  $\Phi$  defined in (\ref{TheCocycle}) by a factor of
\[\exp(H_{\lambda_{1},\lambda_{2}}(w)) = \exp(k(w,\lambda_{1},\lambda_{2})) \exp(il(w,\lambda_{1},\lambda_{2}))
\]
where the functions $k$ and $l$ are defined by 
\begin{equation}\label{eqn:keqn}
k(\lambda_{1},\lambda_{2},\lambda_{3}) = \text{Re}(H_{\lambda_{2},\lambda_{3}}(\lambda_{1}))=  \frac{1}{8}\Big{(}E(\lambda_{1},\lambda_{2},\lambda_{3})+\frac{1}{2}E(i\lambda_{1},i\lambda_{2},\lambda_{3})
+\frac{1}{2}E(i\lambda_{1},\lambda_{2},i\lambda_{3})\Big{)}
\end{equation}
and
\begin{equation}\label{eqn:leqn}
l(\lambda_{1},\lambda_{2},\lambda_{3}) = \text{Im}(H_{\lambda_{2},\lambda_{3}}(\lambda_{1})) = \frac{1}{8}\Big{(}\frac{1}{2}E(\lambda_{1},i\lambda_{2},\lambda_{3})
+\frac{1}{2}E(\lambda_{1},\lambda_{2},i\lambda_{3})- E(i\lambda_{1},\lambda_{2},\lambda_{3})\Big{)}.
\end{equation} 
For future use we record an easy calculation
\begin{equation}\label{eqn:firstfactordef}
l(w,v,i\lambda)-l(w,iv,\lambda) = \frac{1}{8}\Big{(}E(iw,iv,\lambda)-E(iw,v,i\lambda).\Big{)}
\end{equation}
So we have 
\begin{equation}\label{ResultOfTrans}
(w \cdot \Phi)_{\lambda_{1},\lambda_{2}} = \exp(H_{\lambda_{1},\lambda_{2}}(w))\Phi_{\lambda_{1},\lambda_{2}}.
\end{equation}
As is common, in order to put this cocycle in canonical form, we want to kill the part of this extra factor that is multiplication
by a non-zero real number, i.e. multiply it by a boundary to make it unitary.  It is easy to check
that $\lambda_{1} \mapsto \exp(-il(w,v,\lambda_{1}))$ does the job; however it is not holomorphic.  However we
can add another term to make it holomorphic.   Let
$\eta^{w} \in C^{1}(\Lambda, \mathcal{O}^{\times}(V))$ be defined by
\begin{equation}\label{eqn:firstfactor}\eta^{w}_{\lambda}(v)= \exp(-il(w,v,\lambda)-l(w,iv,\lambda)).
\end{equation}

So
\[(\delta \eta^{w})_{\lambda_{1},\lambda_{2}} = \exp(-il(w,\lambda_{1},\lambda_{2}) -l(w,i\lambda_{1},\lambda_{2}))
\]
Therefore
\begin{equation}\label{eqn:firstboundary}
((w \cdot \Phi) (\delta \eta^{w}))_{\lambda_{1},\lambda_{2}}= \Phi_{\lambda_{1},\lambda_{2}} \exp(k(w,\lambda_{1},\lambda_{2}) - l(w,i\lambda_{1},\lambda_{2})).
\end{equation}
Recall that
\[l(w,\lambda_{1},\lambda_{2}) = \frac{1}{8}\Big{(}\frac{1}{2}E(w,i\lambda_{1},\lambda_{2})+\frac{1}{2}E(w,\lambda_{1},i\lambda_{2})- E(iw,\lambda_{1},\lambda_{2})\Big{)}
\]
and
\[k(w,\lambda_{1},\lambda_{2}) =   \frac{1}{8}\Big{(}E(w,\lambda_{1},\lambda_{2})+\frac{1}{2}E(iw,i\lambda_{1},\lambda_{2})
+\frac{1}{2}E(iw,\lambda_{1},i\lambda_{2})\Big{)}
\]
so
\[l(w,i\lambda_{1},\lambda_{2})= \frac{1}{8}\Big{(}-\frac{1}{2}E(w,\lambda_{1},\lambda_{2})+\frac{1}{2}E(w,i\lambda_{1},i\lambda_{2})- E(iw,i\lambda_{1},\lambda_{2})\Big{)}.
\]
Putting this together we can write out the term appearing in (\ref{eqn:firstboundary}) more explicitly:

\begin{equation}
\begin{split}
& k(w,\lambda_{1},\lambda_{2}) - l(w,i\lambda_{1},\lambda_{2}) \\
& = \frac{1}{8}\Big{(}\frac{3}{2}E(w,\lambda_{1},\lambda_{2}) +\frac{3}{2}E(iw,i\lambda_{1},\lambda_{2})
+\frac{1}{2}E(iw,\lambda_{1},i\lambda_{2})-\frac{1}{2}E(w,i\lambda_{1},i\lambda_{2})\Big{)}.
\end{split}
\end{equation}
When we skew-symmetrize $k(w,\lambda_{1},\lambda_{2}) - l(w,i\lambda_{1},\lambda_{2})$, we get the element
 of $Alt^{2}(\Lambda,\mathbb{R})$ given by
\begin{equation}\label{classabove}
\begin{split}
&  2k(w,\lambda_{1},\lambda_{2})-l(w,i\lambda_{1},\lambda_{2})+l(w,i\lambda_{2},\lambda_{1}) \\
&  =\frac{1}{8}\Big{(}3E(w,\lambda_{1},\lambda_{2})+2E(iw,i\lambda_{1},\lambda_{2})-E(w,i\lambda_{1},i\lambda_{2})+2E(iw,\lambda_{1},i\lambda_{2})\Big{)} \\
&  =\frac{1}{8}\Big{(}5E(w,\lambda_{1},\lambda_{2})-3E(w,i\lambda_{1},i\lambda_{2})\Big{)}.
\end{split}
\end{equation}
Therefore, in terms of the data 
\[(B,E) \in \text{Alt}^{2}(\Lambda,\mathbb{R}) \times \text{Alt}^{3}(\Lambda, \mathbb{Z})^{(2,1)+(1,2)}
\] which describes the gerbe (see (\ref{TheCocycle})), the translation looks like
\begin{equation}\label{dataTransform}
\left( B,E \right) \mapsto \left( B+\frac{1}{8}\Big{(}5E(w,\cdot,\cdot)-3E(w,i\cdot,i\cdot)\Big{)}, E \right)
.\end{equation}

Since the element $\frac{1}{8}\Big{(}3E(w,\lambda_{1},\lambda_{2})+3E(w,i\lambda_{1},i\lambda_{2})\Big{)}$ is
clearly of type $(1,1)$ we conclude that the translation of the gerbe is equivalent to the original gerbe tensored with the gerbe described by
\[
(\lambda_{1},\lambda_{2}) \mapsto \exp(E(w,\lambda_{1},\lambda_{2})).
\]

\begin{lem}
For any gerbe $\mathfrak{G}$ on a complex torus $X=V/\Lambda$ defined by a canonical cocycle $\Phi$, with topological type $E$ and any $w \in V$ the canonical cocycle for the gerbe $w \cdot \mathfrak{G}$ is 
\[(\lambda_{1}, \lambda_{2}) \mapsto exp(E(w,\lambda_{1},\lambda_{2})) \Phi_{\lambda_{1},\lambda_{2}}
\]
The translation of $\mathfrak{G}$ by $w$ is isomorphic to
the original if and only if 
\[E(w,\cdot, \cdot)
\in Alt^{2}(\Lambda,\mathbb{Z})+Alt^{2}(\Lambda, \mathbb{R})^{(1,1)}.\]
\end{lem}

\ \hfill $\Box$

For instance, if $w \in \Lambda$ then the translation is actually trivial
as a map of the torus and so this must be the case that this is possible.
In order to avoid writing this element more than necessary, we record this as follows
\begin{dfn}
For a holomorphic gerbe with topological class $E \in
Alt^{3}(\Lambda, \mathbb{Z}) \cong H^{3}(X,\mathbb{Z})$ on a complex torus
$X=V/\Lambda$, we denote by $P$ the map
\[X \to Alt^{2}(\Lambda, \mathbb{R})/(Alt^{2}(\Lambda,\mathbb{Z})+
Alt^{2}(\Lambda,\mathbb{R})^{(1,1)})
\]
given by
\[P(x) = [(\lambda_{1},\lambda_{2}) \mapsto E(w,\lambda_{1},\lambda_{2})]
\]
for $w$ any lift of $x$ to $V$.
\end{dfn}

\begin{dfn} Define the group $K(E,V)$ of a holomorphic gerbe with
topological class $E \in Alt^{3}(\Lambda,\mathbb{Z})$ on a complex
torus $X= V/\Lambda$ to be the subgroup of $V$ defined by
\begin{equation}\label{KEV}K(E,V) = \{w \in V | E(w,\lambda_{1},\lambda_{2}) \in
Alt^{2}(\Lambda,\mathbb{Z})+Alt^{2}(\Lambda,\mathbb{R})^{(1,1)} \}
\end{equation}
or simply
\[K(E,V) = P^{-1}(e_{X}).
\]
\end{dfn}
Thus as a result of the description of gerbes in terms of canonical cocycles from \cite{BGerbes} we have shown the following corollary:
\begin{cor}  Let $w$ be an element in $V$ lifting some translation
action by an element of $X$ on itself.  Then a holomorphic gerbe
with topological class
$E \in Alt^{3}(\Lambda,\mathbb{Z}) \cong H^{3}(X,\mathbb{Z})$ is invariant
under translation by $w$ precisely when the element of
$Alt^{2}(\Lambda, \mathbb{R})$ given by
\[(\lambda_{1},\lambda_{2}) \mapsto E(w,\lambda_{1},\lambda_{2})
\] goes to zero in
\[Alt^{2}(\Lambda, \mathbb{R})/(Alt^{2}(\Lambda,\mathbb{Z})
+Alt^{2}(\Lambda,\mathbb{R})^{(1,1)})\] or equivalently
\[ E(w,\lambda_{1},\lambda_{2}) \in Alt^{2}(\Lambda,\mathbb{Z})+Alt^{2}(\Lambda,\mathbb{R})^{(1,1)}.
\]
\end{cor}
It should be noted that the Hodge Projection 
\[H^{2}(X,\mathbb{R})= Alt^{2}(\Lambda,\mathbb{R}) \to \wedge^{2}(\overline{V}^{\vee})=H^{2}(X,\mathcal{O})\] given by 
\begin{equation}\label{eqn:ActualHodgeProj}
\omega^{H}(w_{1}, w_{2}) = \frac{1}{4}\bigg{(}\omega(w_{1}, w_{2})-\omega(iw_{1}, iw_{2})+i\omega(iw_{1}, w_{2})+i\omega(w_{1}, iw_{2})\bigg{)}
\end{equation}
sets up an isomorphism 
\[Alt^{2}(\Lambda, \mathbb{R})/(Alt^{2}(\Lambda,\mathbb{Z}) +Alt^{2}(\Lambda,\mathbb{R})^{(1,1)})\cong \wedge^{2}\overline{V}^{\vee}/Alt^{2}(\Lambda,\mathbb{Z})^{H}
\]
where $Alt^{2}(\Lambda,\mathbb{Z})^{H}$ is the image of $Alt^{2}(\Lambda,\mathbb{Z})$ under the Hodge projection.

Notice that $\Lambda \subset K(E,V)$ as expected from the interpretation
of $K(E,V)$ as automorphisms fixing a gerbe.  Indeed for $w=\lambda_{3} \in \Lambda$
we can add an element of type $(1,1)$ to make the term appearing in Equation \ref{classabove}  integral:
\begin{equation}
\begin{split}
& \frac{1}{8}\Big{(}5E(\lambda_{3},\lambda_{1},\lambda_{2})-3E(\lambda_{3},i\lambda_{1},i\lambda_{2})\Big{)} +
\frac{1}{8}\Big{(}3E(\lambda_{3},\lambda_{1},\lambda_{2})+3E(\lambda_{3},i\lambda_{1},i\lambda_{2})\Big{)} \\
& = E(\lambda_{3},\lambda_{1},\lambda_{2}) \in \mathbb{Z}.
\end{split}
\end{equation}
In fact more generally we have
\[E(w,\lambda_{1},\lambda_{2}) \in \mathbb{Z} \text{  } \forall  \lambda_{1},\lambda_{2} \in \Lambda \Longrightarrow w \in K(E,V).
\]
and therefore $\Lambda \subset K(E,V)$.  Thus we have shown the following
\begin{lem}\label{explicit}
Given any complex torus $X= V/\Lambda$ and a holomorphic gerbe $\mathfrak{G}$ on $X,$
\[K(\mathfrak{G}) = K(E,V)/\Lambda \subset X
\]
where $K(\mathfrak{G})$ was defined in definition (\ref{Kdef}) and $K(E,V)$ was defined in Equation (\ref{KEV}).
\end{lem}
\begin{ex}\label{finiteHeis}
Let $X= V/\Lambda$ be a complex torus of complex dimension $2$ such that the Neron-Severi group of $X$ is trivial.  In this case $P(w) = E(w,\cdot, \cdot)$.  Choose a basis $e_{1}, e_{2}, e_{3}, e_{4}$ of $V$ and let 
\[E =  e_{1}^{*}\wedge e_{2}^{*} \wedge e_{3}^{*}.
\] 
We compute here 
\[\{w|E(w,\cdot, \cdot) \in \text{Alt}^{2}(\Lambda, \mathbb{Z}) \}/\Lambda \subset K(\mathfrak{G}).
\]
where $\mathfrak{G}$ is a gerbe of topological type $E$.
Given an element  
\[w = w_{1}e_{1}+ w_{2} e_{2} + w_{3} e_{3} + w_{4} e_{4}  \in (\mathbb{R}/\mathbb{Z})^{4},
\]
$E(w,\cdot, \cdot)$  lives in $\text{Alt}^{2}(\Lambda, \mathbb{Z})$ if and only if $w_1,w_2,w_3 \in \mathbb{Z}$ and therefore $S^{1}\subset K(\mathfrak{G})$ and $\mathbb{R} \subset \Lambda(\mathfrak{G})$.  For another example we can take $E$ to be twice any class in $Alt^{3}(\Lambda,\mathbb{Z})$.  Let $\mathfrak{G}$ be a holomorphic gerbe on $X$ with topological type $E$.  Then we have a finite subgroup $\frac{1}{2}\Lambda/\Lambda \subset K(\mathfrak{G})$ and a discrete subgroup $\frac{1}{2}\Lambda \subset \Lambda(\mathfrak{G})$. 
 \end{ex}

\section{An Explicit Choice of Isomorphism}\label{ExplicitChoice}

In the following we would like to actually pick an equivalence for each element of $\Lambda(\mathfrak{G})$.  In practical terms, given $w \in \Lambda(\mathfrak{G})$, we would like to find $\tau^{w} \in C^{1}(\Lambda, \mathcal{O}^{\times}(V))$ such that
\[(w \cdot \psi) (\delta \tau^{w})= \psi.
\]

\begin{dfn}
For every translation preserving a gerbe (presented by a canonical cocycle) we now choose an explicit isomorphism from the translated gerbe to the original.  In order to do this we chose a decomposition of the term appearing in (\ref{dataTransform}).
\begin{equation} \label{eqn:decomp}  \frac{1}{8}\Big{(}5E(w,\lambda_{1},\lambda_{2})-3E(w,i\lambda_{1},i\lambda_{2})\Big{)} = E^{w}(\lambda_{1},\lambda_{2}) +\epsilon^{w}(\lambda_{1},\lambda_{2})
\end{equation}
where $E^{w} \in Alt^{2}(\Lambda,\mathbb{R})^{(1,1)}$ and $\epsilon^{w} \in Alt^{2}(\Lambda,\mathbb{Z})$,  corresponding to a gerbe that is translation invariant.   This decomposition is unique when the Neron-Severi group $ Alt^{2}(\Lambda,\mathbb{R})^{(1,1)} \cap  Alt^{2}(\Lambda,\mathbb{Z})$ of the torus is zero, but in general it is not unique.  The decomposition in Equation (\ref{eqn:decomp}) is equivalent to the two conditions 
\[E(w, \cdot, \cdot)^{(2,0)+(0,2)} = (\epsilon^{w})^{(2,0)+(0,2)}
\]
and 
\[\frac{1}{4}E(w,\cdot, \cdot)^{(1,1)} = E^{w}+ (\epsilon^{w})^{(1,1)}
\]
\end{dfn}

Define $\mu^{w} \in C^{1}(\Lambda,U(1))$ by
\begin{equation}\label{eqn:secondfactor}
\mu^{w}_{\lambda}= \exp \left(\frac{-1}{16}\Big{(}\frac{3}{2}E(iw,i\lambda,\lambda)+\frac{1}{2}E(iw,\lambda,i\lambda)\Big{)} \right).
\end{equation}
We can use $\mu^{w}$ as a first step towards bounding the term multiplying $\Phi$ in Equation (\ref{eqn:firstboundary}):

\begin{equation}\label{eqn:secondboundary}
\begin{split}
 & \exp\Big{(}k(w,\lambda_{1},\lambda_{2})-l(w,i\lambda_{1},\lambda_{2})\Big{)} (\delta \mu^{w})_{\lambda_{1},\lambda_{2}} \\
& =
 \exp \left(\frac{1}{16}\Big{(}5E(w,\lambda_{1},\lambda_{2})-3E(w,i\lambda_{1},i\lambda_{2})\Big{)} \right) \\
& = \exp \left(\frac{1}{2}E^{w}(\lambda_{1},\lambda_{2}) + \frac{1}{2}\epsilon^{w}(\lambda_{1},\lambda_{2}) \right).
\end{split}
\end{equation}

We know that this is a boundary precisely because its skew symmetrization lives in $Alt^{2}(\Lambda,\mathbb{Z})+Alt^{2}(\Lambda,\mathbb{R})^{(1,1)}$.
We now produce elements $\nu^{w}$ and  $\phi^{w}$ in $C^{1}(\Lambda,\mathcal{O}^{\times}(V))$, bounding the two terms $\exp \left(\frac{1}{2}\epsilon^{w}(\lambda_{1},\lambda_{2}) \right)$ and $\exp\left(\frac{1}{2}E^{w}(\lambda_{1},\lambda_{2})\right)$ which appear on the right hand side of  Equation (\ref{eqn:secondboundary}).  If
\[\lambda = \sum_{i} n_{i} \lambda^{i}
\]
\begin{equation} \label{eqn:thirdfactor}
\nu^{w}_{\lambda} = \exp\left(-\frac{1}{2}\sum_{i<j}\epsilon^{w}(n_{i} \lambda^{i}, n_{j} \lambda^{j}) \right)
\end{equation}
then
\begin{equation}\label{eqn:thirdboundary}
\exp \left(\frac{1}{2}\epsilon^{w}(\lambda_{1},\lambda_{2}) \right) (\delta \nu)_{\lambda_{1},\lambda_{2}} = 1.
\end{equation}

Let $\phi^{w} \in C^{1}(\Lambda,\mathcal{O}^{\times}(V))$ be defined by
\begin{equation} \label{eqn:fourthfactor}
\phi^{w}_{\lambda} = \exp(L^{w}(v,\lambda))
\end{equation}
where
\[L^{w}(v,\lambda) = \frac{i}{2}E^{w}(iv,\lambda)-\frac{1}{2}E^{w}(v,\lambda)+\frac{i}{4}E^{w}(i\lambda,\lambda ).
\]
Then
\begin{equation}\label{eqn:fourthboundary}
\exp\left(\frac{1}{2}E^{w}(\lambda_{1},\lambda_{2})\right) (\delta \phi^{w})_{\lambda_{1},\lambda_{2}} = 1,
\end{equation}
here the choice of $L^{w}$ can easily be guessed from the Appell-Humbert theorem for line bundles discussed in \cite{BGerbes}.

Combining Equations (\ref{eqn:firstfactor}), (\ref{eqn:secondfactor}), (\ref{eqn:thirdfactor}), and (\ref{eqn:fourthfactor}) we define
\begin{equation}\label{eqn:tauDef}
\tau^{w}_{\lambda}= \phi^{w}_{\lambda} \nu^{w}_{\lambda} \mu^{w}_{\lambda} \eta^{w}_{\lambda}
\end{equation}
which gives an isomorphism
\[w \cdot \mathfrak{G} \to \mathfrak{G}
\]
 from the translation of the gerbe back to the original gerbe.  The formula giving $\tau$ is 
\begin{equation}\label{tauequation}\tau^{w}_{\lambda}(v) = \exp(T(w,v,\lambda))
\end{equation} where
\begin{equation}
\begin{split}
T(w,v,\lambda) & =   L^{w}(v,\lambda) -\frac{1}{16}\left(\frac{3}{2}E(iw,i\lambda,\lambda)+\frac{1}{2}E(iw,\lambda,i\lambda)\right) \\
&  \quad   - il(w,v,\lambda)-l(w,iv,\lambda) -\frac{1}{2}\sum_{i<j}\epsilon^{w}(n_{i} \lambda^{i}, n_{j} \lambda^{j}).
\end{split}
\end{equation}

Using Equations (\ref{eqn:firstboundary}), (\ref{eqn:secondboundary}), (\ref{eqn:thirdboundary}), and (\ref{eqn:fourthboundary}) we conclude that for $w \in K(E)$ we have
\[\exp(H_{\lambda_{1},\lambda_{2}}(w)) (\delta \tau^{w})_{\lambda_{1},\lambda_{2}} = 1
\]
and so from the behavior under translation given in Equation (\ref{groupext2}) we have
\[((w \cdot \Phi) (\delta \tau^{w}))_{\lambda_{1},\lambda_{2}}= \Phi_{\lambda_{1},\lambda_{2}}.
\]

\section{The First Obstruction}\label{FirstObs}

The map from $\mathfrak{G}$ to itself shown in (\ref{GerbeSequence})  is
the $\mathcal{O}^{\times}$-torsor $\mathcal{L}^{w_{1},w_{2}}$ on $X$ described by the cocycle $\Xi^{w_{1},w_{2}}$ in $Z^{1}(\Lambda,\mathcal{O}^{\times}(V))$
\[\Xi^{w_{1},w_{2}}: \Lambda \to \mathcal{O}^{\times}(V)
\]
given by

\begin{equation}\label{xidef} \Xi^{w_{1},w_{2}}= (\tau^{w_{1}+w_{2}})^{-1}(w_{1} \cdot \tau^{w_{2}})\tau^{w_{1}}.
\end{equation}

In order to proceed, recall that we choose for each $w \in \Lambda(\mathfrak{G})$ a decomposition as in Equation (\ref{eqn:decomp}).   Fix $w_{1}$ and $w_{2}$ and write
\[E^{w_{1}+w_{2}} = E^{w_{1}} + E^{w_{2}} + C^{w_{1},w_{2}}
\]
and
\[\epsilon^{w_{1}+w_{2}} = \epsilon^{w_{1}} + \epsilon^{w_{2}} +\zeta^{w_{1},w_{2}}\]
where $C=- \zeta \in NS(X)$.
Notice that $NS(X)$ is acted on trivially by translations and $C$ is also symmetric under exchange of $w_{1}$ and $w_{2}$.  
\begin{ass}\label{ass}
We can easily find several cases where it is clear that the $E^{w}$ can be chosen so that $C$ vanishes (it could be that this is always the case).   We will from now only consider subgroups of $\Lambda(\mathfrak{G})$
for which $C$ vanishes.  \end{ass}
Notice that above assumption holds in the following cases

\begin{enumerate}
\item when we restrict to the subgroup $\Lambda(\mathfrak{G})_{\mathbb{Z}}$ of $\Lambda(\mathfrak{G})$ 
defined by \[\Lambda(\mathfrak{G})_{\mathbb{Z}}=\{w|E(w,\cdot,\cdot) \in Alt^{2}(\Lambda,\mathbb{Z}) \},\]
this will be called the integral case,
\item when we restrict  to the subgroup of $\Lambda(\mathfrak{G})_{(1,1)}$ of $\Lambda(\mathfrak{G})$ defined by \[\Lambda(\mathfrak{G})_{(1,1)} = \{w|E(w,\cdot, \cdot) \in Alt^{2}(\Lambda, \mathbb{R})^{(1,1)} \},\]  
this will be called the case of type $(1,1)$,
\item in the case where $NS(X)$ vanishes and we work 
with the full group $\Lambda(\mathfrak{G})$.
\end{enumerate}

In the remainder of the paper, we will calculate the classes $[\Xi]$ in  $H^{2}(\Lambda(\mathfrak{G})_{\mathbb{Z}}, \text{Pic}(X))$  

\noindent
and  $H^{2}(\Lambda(\mathfrak{G})_{(1,1)}, \text{Pic}(X))$  corresponding to the first two of the above cases.  This class is the first obstruction to equivariance.  When this class  vanishes, we calculate a second (and final) obstruction to equivariance, also in the first two of the above cases.  The second obstructions in the first two of the above cases will live in $H^{3}(\Lambda(\mathfrak{G})_{\mathbb{Z}},\mathbb{C}^{\times})$ and $H^{3}(\Lambda(\mathfrak{G})_{(1,1)},\mathbb{C}^{\times})$ respectively.  The first and second obstructions correspond to the images of $[\mathfrak{G}]$ under the second and third differentials from the spectral sequence (\ref{ss}).

All terms contributing to $\tau$  (see Equation (\ref{eqn:tauDef})) which do not
depend on $v$ make no contribution to the
factor $\Xi^{w_{1},w_{2}}$.  This makes $\mu$, $\nu$ and part of
$\phi$ irrelevant and so we are left with
$\eta$ (see Equation (\ref{eqn:firstfactor})) and the non-constant part
of $\phi$ (see Equation (\ref{eqn:fourthfactor})).

\begin{equation}\label{XOBS}
\begin{split}
 \Xi^{w_{1},w_{2}}_{\lambda}
& = (\tau^{w_{1}+w_{2}}_{\lambda})^{-1}(w_{1} \cdot \tau^{w_{2}}_{\lambda})\tau^{w_{1}}_{\lambda} \\
&=\exp\left( T(w_{2},v+w_{1},\lambda) - T(w_{1}+w_{2},v,\lambda) 
+T(w_{1},v,\lambda)\right)   \\
& = \exp\left(-L^{w_{1}+w_{2}}_{\lambda} +w_{1} \cdot L^{w_{2}}_{\lambda} + L^{w_{1}}_{\lambda} -il(w_{2},w_{1},\lambda)-l(w_{2},iw_{1},\lambda)\right) \\
& = \exp\left(\frac{i}{2}E^{w_{2}}(iw_{1},\lambda)-\frac{1}{2}E^{w_{2}}(w_{1},\lambda)
-il(w_{2},w_{1},\lambda)-l(w_{2},iw_{1},\lambda)\right) \\
& =  \exp(S(w_{1},w_{2},\lambda))
\end{split}
\end{equation}

where

\begin{equation} \label{Seqn}
\begin{split}
S(w_{1},w_{2},\lambda)= &\frac{i}{2}E^{w_{2}}(iw_{1},\lambda)-\frac{1}{2}E^{w_{2}}(w_{1},\lambda)
-il(w_{2},w_{1},\lambda)-l(w_{2},iw_{1},\lambda) \\
= & -\frac{1}{2}E^{w_{2}}(w_{1},\lambda) -  \frac{1}{8}\Big{(}-\frac{1}{2}E(w_{2},w_{1},\lambda)
+\frac{1}{2}E(w_{2},iw_{1},i\lambda)- E(iw_{2},iw_{1},\lambda)\Big{)} \\
& + \frac{i}{2}E^{w_{2}}(iw_{1},\lambda) -  \frac{i}{8}\Big{(}\frac{1}{2}E(w_{2},iw_{1},\lambda)
+\frac{1}{2}E(w_{2},w_{1},i\lambda)- E(iw_{2},w_{1},\lambda)\Big{)}
\end{split}
\end{equation}

Recall that $l$ was defined in Equation (\ref{eqn:leqn}) and for $w \in \Lambda(\mathfrak{G})$, $E^{w}$ is given by Equation (\ref{eqn:decomp}).  Notice that $\Xi$ is constant, it does not depend on $v$.  It is also easy to see that as a cohomology class in $H^{2}(\Lambda(\mathfrak{G}), H^{1}(X,\mathcal{O}^{\times}))$ that $\Xi$ does not depend on the choices of $E^{w}$.  As an element in $Z^{1}(\Lambda,\mathcal{O}^{\times}(V))$, $\Xi^{w_{1},w_{2}}$
is equivalent to an element of
$Z^{1}(\Lambda,U(1)) = \text{Hom}(\Lambda,U(1))$.    Therefore $[\Xi]$ is in the image of the obvious map $H^{2}(\Lambda(\mathfrak{G}), \text{Pic}^{0}(X)) \to H^{2}(\Lambda(\mathfrak{G}), \text{Pic}(X))$.

\begin{rmk}
The multiplication rule in the group $\pi_{0}(G(\mathfrak{G}))$ described in (\ref{LatticeSES}) is 
\begin{equation}
\begin{split}
& \Big{(}\alpha_{1},w_{1}\Big{)}\Big{(}\alpha_{2},w_{2}\Big{)}  = \\
& \Big{(}\alpha_{1} \alpha_{2}(\tau^{w_{1}+w_{2}}_{\lambda})^{-1}(w_{1}
\cdot \tau^{w_{2}}_{\lambda})\tau^{w_{1}}_{\lambda} , x_{1}+x_{2} \Big{)} =
\bigg{(}\alpha_{1}\alpha_{2} \exp\Big{(}S(w_{1},w_{2},\cdot)\Big{)},w_{1}+w_{2}\bigg{)}
\end{split}
\end{equation}
Its easy to see directly that this is associative, but this also follows
from the functoriality of pulling back gerbes.
\end{rmk}

When we multiply (\ref{XOBS}) by the boundary of $\Theta \in C^{2}(K(\mathfrak{G}),\mathcal{O}^{\times}(V))$ defined by 
\begin{equation}\label{ThetaEqn}\Theta^{w_{1},w_{2}}(v) = \exp \left(il(w_{2},w_{1},v)+l(w_{2},w_{1},iv)
-\frac{i}{2}E^{w_{2}}(iw_{1},v)-\frac{1}{2}E^{w_{2}}(iw_{1},iv) \right)
\end{equation} we bring
it into $Z^{1}(\Lambda,U(1))$ where it becomes
(using Equation (\ref{eqn:firstfactordef}))
\begin{equation} \label{this}
\begin{split}
&  \exp\Big{(}-l(w_{2},iw_{1},\lambda)+l(w_{2},w_{1},i\lambda)
-E^{w_{2}}(w_{1},\lambda)\Big{)} \\
&  = \exp \left(\frac{1}{8}\Big{(}E(iw_{2},iw_{1},\lambda)
-E(iw_{2},w_{1},i\lambda)\Big{)}-E^{w_{2}}(w_{1},\lambda) \right).
\end{split}
\end{equation} 

\begin{cor}
The short exact sequence (\ref{LatticeSES}) is classified by the element of $H^{2}(\Lambda(\mathfrak{G}),\text{Pic}(X))$
represented by
\[(w_{1}, w_{2}) \mapsto \Big{[}\lambda \mapsto \exp \left(\frac{1}{8}\Big{(}E(iw_{2},iw_{1},\lambda)
-E(iw_{2},w_{1},i\lambda)\Big{)}-E^{w_{2}}(w_{1},\lambda) \right)\Big{]}.
\]
This is the first of two possible obstructions towards the equivariance of $\mathfrak{G}$ under $K(\mathfrak{G})$.  When this obstruction vanishes there is an additional obstruction in 
$H^{3}(K(\mathfrak{G}), \mathbb{C}^{\times})$ which we discuss later.
\end{cor}

We now analyze our lift of $[\Xi]$ to $H^{2}(\Lambda(\mathfrak{G}), \text{Pic}^{0}(X))$, where we identify $\text{Pic}^{0}(X) = Hom(\Lambda, U(1))$.
Skew-symmetrizing Equation (\ref{this})
gives us the element in $\text{Alt}^{2}(\Lambda(\mathfrak{G}), Hom(\Lambda, U(1)))$
 given by

\begin{equation}\label{groupext1}
\begin{split}
 & \exp\left(  \frac{1}{8}\Big{(}2E(iw_{2},iw_{1},\lambda)-E(iw_{2},w_{1},i\lambda)
-E(w_{2},iw_{1},i\lambda)\Big{)}+(-E^{w_{2}}(w_{1},\lambda)+E^{w_{1}}(w_{2},\lambda)) \right) \\
 & = \exp \left( \frac{1}{8}\Big{(}2E(iw_{2},iw_{1},\lambda)-E(w_{2},w_{1},\lambda)
+E(iw_{2},iw_{1},\lambda)\Big{)}+(-E^{w_{2}}(w_{1},\lambda)+E^{w_{1}}(w_{2},\lambda)) \right) \\
 & = \exp \left(  \frac{1}{8}\Big{(}3E(iw_{2},iw_{1},\lambda) - E(w_{2},w_{1},\lambda)\Big{)} +(-E^{w_{2}}(w_{1},\lambda)+E^{w_{1}}(w_{2},\lambda)) \right)
\end{split}
\end{equation}

In the next two sub-sections we refer to the two different cases discussed in assumption \ref{ass}.
\subsection{The Integral Case}
In the first case we look at the subgroup $\Lambda(\mathfrak{G})_{\mathbb{Z}}$ of $w \in \Lambda(\mathfrak{G})$ such that the contraction with $E$ is integral.  In this case we can take 
\begin{equation}\label{fristcase}
E^{w}= -\frac{3}{8}(E(w,\cdot,\cdot)+E(w,i\cdot,i\cdot))\end{equation} and $\epsilon^{w} = E(w,\cdot,\cdot)$.  These are clearly linear in $w$.
In the special case that $w_{1}$ and $w_{2}$ both contract with $E$ to form something integral we get using Equation (\ref{fristcase}) that Equation (\ref{groupext1}) becomes 

\begin{equation}\label{BothInt}
\begin{split}
&\exp \left( \frac{1}{8}\Big{(}3E(iw_{2},iw_{1},\lambda) - E(w_{2},w_{1},\lambda) \Big{)}  \right) \\ & \exp \left ( \frac{3}{8} \Big{(}(E(w_{2},w_{1},\lambda) +E(w_{2},iw_{1},i\lambda) -E(w_{1},w_{2},\lambda) - E(w_{1},iw_{2},i\lambda))\Big{)}  \right) \\
&=\exp \left( \frac{1}{8}\Big{(}3E(iw_{2},iw_{1},\lambda) + 5E(w_{2},w_{1},\lambda)  +3(E(w_{2},iw_{1},i\lambda)  - E(w_{1},iw_{2},i\lambda))\Big{)}  \right) \\
&=\exp \left( E(w_{2},w_{1},\lambda)  \right) 
\end{split}
\end{equation}
Therefore the obstruction in this case is the element of $\text{Alt}^{2}(\Lambda(\mathfrak{G})_{\mathbb{Z}}, \text{Pic}^{0}(X))$ given by 
\[(w_{1}, w_{2}) \mapsto [\lambda \mapsto E(w_{2},w_{1},\lambda)].
\]
We now give an example of a subgroup on which the theta extension (\ref{fristextobs}) of $K(\mathfrak{G})$ by $\text{Pic}(X)$ does not split and therefore the gerbe is not equivariant for this subgroup.
\begin{ex}
Let $X=V/\Lambda$ be any complex torus of dimension $2$.  Let $\Lambda \subset V$ be spanned by $e_{1}, e_{2}, e_{3}$ and $e_{4}$ and let 
\[E = 2e_{1}^{*} \wedge e_{2}^{*} \wedge e_{3}^{*}
\]
Then the contractions of $E$ by $\frac{1}{2}e_{1}$ and $\frac{1}{2}e_{2}$ are clearly integral.   Let $\mathfrak{G}$ be a gerbe on $X$ with topological type $E$.  The vectors $\frac{1}{2}e_{1}$ and $\frac{1}{2}e_{2}$ generate a subgroup $\mathbb{Z}/2\mathbb{Z} \times \mathbb{Z}/2\mathbb{Z} \subset K(\mathfrak{G})$ and it pulls back to a rank $2$ sub-lattice of $\frac{1}{2}\Lambda \subset V$.  Over this sub-lattice the sequence (\ref{LatticeSES}) does not split and so $\mathfrak{G}$ is not equivariant for this copy of 
 $\mathbb{Z}/2\mathbb{Z} \times \mathbb{Z}/2\mathbb{Z}$.  This follows immediately from Remark \ref{annoying} and the fact that (\ref{BothInt}) is not trivial if we evaluate on $\lambda= e_{3}$, $w_{1}= \frac{1}{2}e_{1}$, $w_{2}=\frac{1}{2}e_{2}$.
\end{ex} 
\subsection{The Case Of Type $(1,1)$}
In the second case found in Assumption \ref{ass}, the contraction of $E$ with with $w$ is assumed to be of type $(1,1)$.   We can define 
\[E^{w} = 4E(w,\cdot,\cdot)\] again ensuring that $C$ and $\zeta$ vanish.
Now using Equations (\ref{str1}) and (\ref{eqn:decomp}) we have in the special case that $\epsilon^{w}$ vanishes meaning that $E(w,\cdot, \cdot)$ has type $(1,1)$.
\begin{equation}\label{groupext2}
\begin{split}
& -E^{w_{2}}(w_{1},\lambda)+E^{w_{1}}(w_{2},\lambda) \\
& = \frac{1}{8}\Big{(}-5E(w_{2},w_{1},\lambda) +   3E(w_{2},iw_{1},i\lambda) + 5E(w_{1},w_{2},\lambda) - 3E(w_{1},iw_{2},i\lambda)\Big{)} \\
& = \frac{1}{8}\Big{(}-10E(w_{2},w_{1},\lambda) + 3E(w_{2},w_{1},\lambda) -3E(iw_{2},iw_{1},\lambda)\Big{)} \\
& = \frac{1}{8}\Big{(}-7E(w_{2},w_{1},\lambda) -3E(iw_{2},iw_{1},\lambda)\Big{)}
\end{split}
\end{equation}
Therefore substituting Equation (\ref{groupext2}) into Equation (\ref{groupext1}) we see the factor is
\begin{equation}
\begin{split}
 & \exp \left( \frac{1}{8}\Big{(}-8E(w_{2},w_{1},\lambda) \Big{)}\right) \\
 & = \exp\Big{(}E(w_{1},w_{2},\lambda)\Big{)}.
 \end{split}
\end{equation}
Therefore the second obstruction in this case is the element of $\text{Alt}^{2}(\Lambda(\mathfrak{G})_{(1,1)}, \text{Pic}^{0}(X))$ given by 
\[(w_{1}, w_{2}) \mapsto [\lambda \mapsto E(w_{1},w_{2},\lambda)].
\]

\section{The Second Obstruction}\label{SecondObs}
When the first obstruction vanishes then the element $\Theta$ defined in Equation (\ref{ThetaEqn}) can be considered a trivialization of the first obstruction (whenever such a trivialization exists).  

\begin{rmk}Strictly speaking to fully trivialize the first obstruction $\Xi$ one must also choose a cocycle whose boundary is the difference between (\ref{eqn:firstfactordef}) and its skew symmetrization.  However, this comes at no cost because such a cocycle will have no boundary with respect to the group $\Lambda(\mathfrak{G})$.
\end{rmk}
The boundary map of $\Theta$ with respect to the group $\Lambda(\mathfrak{G})$ gives the following cocycle in $Z^{3}(\Lambda(\mathfrak{G}),\mathbb{C}^{\times})$: 
\begin{equation}\label{EvalTheta}
\Theta^{w_{2}, w_{3}}(w_{1}) = \exp \left(il(w_{3},w_{2},w_{1})+l(w_{3},w_{2},iw_{1})
-\frac{i}{2}E^{w_{3}}(iw_{2},w_{1})-\frac{1}{2}E^{w_{3}}(iw_{2},iw_{1}) \right)
\end{equation}
The associated class in $H^{3}(\Lambda(\mathfrak{G}),\mathbb{C}^{\times})$ represents the lift to $\Lambda(\mathfrak{G})$ of the image of $[\mathfrak{G}]$ under the differential $d_{3}$ from the spectral sequence introduced in (\ref{ss}).
When we skew-symmetrize this with respect to $w_{1}, w_{2}, w_{3}$ we get the image of the second obstruction under the isomorphism 
\[H^{3}(\Lambda(\mathfrak{G}),\mathbb{C}^{\times}) \cong \text{Alt}^{3}(\Lambda(\mathfrak{G}), \mathbb{C}^{\times}).\]  It comes out to be 
\begin{equation} \label{GeneralFactor}
\exp \left(3 \left(E^{w_{3}}(w_{1},w_{2}) + E^{w_{1}}(w_{2},w_{3})-E^{w_{2}}(w_{1},w_{3}) \right)\right)
\end{equation}

\subsection{The Integral Case}
In the integral case we restrict to the subgroup $\Lambda(\mathfrak{G})_{\mathbb{Z}}$ defined in Assumption \ref{ass} and we make the substitution 

\[E^{w}= -\frac{3}{8}(E(w,\cdot,\cdot)+E(w,i\cdot,i\cdot))\]
the expression (\ref{GeneralFactor}) reduces to 
\begin{equation}\label{nine}[(w_{1},w_{2},w_{3}) \mapsto \exp(-9E(w_{1}, w_{2},w_{3})) ] \in \text{Alt}(\Lambda(\mathfrak{G})_{\mathbb{Z}},\mathbb{C}^{\times})
\end{equation}
\begin{ex}
Let $X=V/\Lambda$ be two dimensional with $NS(X)=0$ and a basis for $V$ given by $e_{1},e_{2},e_{3},e_{4}$, let $\mathfrak{G}$ be a gerbe with topological type 
\[E = 4 e_{1}^{*} \wedge e_{2}^{*} \wedge e_{3}^{*}.
\]
If we consider the subgroup $\frac{1}{2} \Lambda \subset \Lambda(\mathfrak{G})$ then the first obstruction to $\mathfrak{G}$ being equivariant for this subgroup vanishes.  The assumption that $NS(X)=0$ ensures using the Hochschild-Serre spectral sequence that the map 
\[H^{2}(\frac{1}{2} \Lambda/\Lambda, \text{Pic}(X)) \to H^{2}(\frac{1}{2} \Lambda, \text{Pic}(X))
\]
is injective and so the first obstruction to $\frac{1}{2} \Lambda/\Lambda$-equivariance vanishes as well.  Applying Remark \ref{annoying} and Equation (\ref{nine}) we see that the second obstruction for $\frac{1}{2} \Lambda$-equivariance does not vanish because $8$ does not divide $4\cdot9=36$.  Therefore the second obstruction for $\frac{1}{2} \Lambda/\Lambda$-equivariance also does not vanish.
\end{ex}

\subsection{The Case of Type $(1,1)$}
In this case we restrict to the subgroup $\Lambda(\mathfrak{G})_{(1,1)}$ defined in Assumption \ref{ass} and we can make the substitution \[E^{w} = 4E(w,\cdot,\cdot)\]  into (\ref{GeneralFactor}) to reduce it to
\[[(w_{1},w_{2},w_{3}) \mapsto \exp(36E(w_{1},w_{2},w_{3}))] \in \text{Alt}(\Lambda(\mathfrak{G})_{(1,1)},\mathbb{C}^{\times}).
\]

\section{A Gerbal Action on the Stack of Twisted Sheaves}\label{Gerbal Action}
 There sheaves of $\mathcal{O}_{\mathfrak{G}}$-modules on a gerbe $\mathfrak{G}$ form an abelian category $\text{Mod}(\mathcal{O}_{\mathfrak{G}})$.  We would like to think of this category as a sort of categorical analogue of a vector space.  The direct sum here plays the role of addition while the action of the groupoid $\mathcal{P}ic^{0}(X)$ plays the role of scalar multiplication.  The group $K(\mathfrak{G})$ defined abstractly in Definition \ref{Kdef} and explicitly described in Lemma \ref{explicit} of translations of a complex torus $X$ that fix a gerbe  $\mathfrak{G}$ on  $X$ act on this category by functors which preserve the structure of direct sum and tensor by elements of $\mathcal{P}ic^{0}(X)$.  There are three possible ways that it can act which we describe in order from weakest to strongest.  Each case is a special case of the ones preceding it.   Following Frenkel and Zhu from \cite{FZ} we call the middle case a gerbal representation.  These actions are analogous to the projective representation on $H^{0}(X,\mathcal{L})$ (or twisted equivariant structure on $\mathcal{L}$) of the group of translations which fix a holomorphic line bundle $\mathcal{L}$ on a complex torus $X$.  Although we do not try it here, it might be interesting to look at twisted equivariant sheaves in this context instead of just twisted sheaves.

Let $\mathfrak{G}$ be an $\mathcal{O}^{\times}$ gerbe on a complex torus $X=V/\Lambda$ with topological class 
\[E \in  \text{Alt}^{3}(\Lambda, \mathbb{Z}) = H^{3}(X,\mathbb{Z}).\]   
Let $D \subset K(\mathfrak{G})$ be a discrete subgroup following 
the requirements in Assumption \ref{ass}.
Then for every $x \in D$ there is an autoequivalence $F_{x}$ of $\text{Mod}(\mathcal{O}_{\mathfrak{G}})$ as an abelian category which commutes with the action of $\mathcal{P}ic^{0}(X)$.   For every pair of elements $x_{1},x_{2} \in D$ there are invertible natural transformations 
\[N_{x_{1},x_{2}}:F_{x_{1}+x_{2}} \to T_{\mathcal{L}^{x_{1},x_{2}}} \circ F_{x_{1}} \circ F_{x_{2}}
\]
where $T_{\mathcal{L}^{x_{1},x_{2}}}$ is the functor given by tensoring with the line bundles corresponding to the $\mathcal{L}^{x_{1},x_{2}}$.  These line bundles are analogous to the discrepancy that appears in the failure of a projective representation to be an actual representation.   In the standard case of line bundles on a torus, the fact that it is a scalar is due to compactness of the torus but in general situations it need not be so.  In our case, the fact that the line bundles $\mathcal{L}$ are of degree zero follows from Assumption \ref{ass}, but in general this should be generalized to arbitrary line bundles.   It easy to see in this case that for every triple $x_{1}, x_{2}, x_{3}$ there is an isomorphism 
\[\mathcal{L}^{x_{1}, x_{2}+x_{3}} \otimes
\mathcal{L}^{x_{2}, x_{3}}
\cong \mathcal{L}^{x_{1}, x_{2}} \otimes \mathcal{L}^{x_{1}+x_{2}, x_{3}}
\]
which satisfies a coherence condition involving $4$ group elements.   
This means that $\mathcal{L}$ is in the category of $2$-cocycles of the group $D$ with values in the groupoid $\mathcal{P}ic^{0}(X)$ as defined in \cite{Po2008}.

 For the next strongest notion of group action we assume that the line bundles can be trivialized.  This can be done under the condition of vanishing of the first obstruction.  Therefore we consider invertible natural transformations 
\[N_{x_{1},x_{2}}:F_{x_{1}+x_{2}} \to  F_{x_{1}} \circ F_{x_{2}}
\]

Given $x_{1}, x_{2}, x_{3} \in D$ there are two ways of mapping  $F_{x_{1}+x_{2}+x_{3}}$ to $F_{x_{1}} \circ F_{x_{2}} \circ F_{x_{3}}$:

\begin{equation}
\label{important}
\xymatrix{**[l]{F_{x_{1}+x_{2}+x_{3}}}  \ar@/_1.5pc/[r]_-{
(N_{x_{1},x_{2}} \circ \text{id}_{F_{x_{3}}})\circ N_{x_{1}+x_{2},x_{3}}
} \ar@/^1.5pc/[r]^-{
(\text{id}_{F_{x_{1}}} \circ N_{x_{2},x_{3}})\circ N_{x_{1},x_{2}+x_{3}}
} & **[r]  F_{x_{1}} \circ F_{x_{2}} \circ F_{x_{3}} }
\end{equation}

The two natural transformations shown above factor as  
\[F_{x_{1}+x_{2}+x_{3}} \stackrel{N_{x_{1}+x_{2},x_{3}}}\longrightarrow F_{x_{1}}  \circ F_{x_{2}+ x_{3}} 
\stackrel{N_{x_{1},x_{2}} \circ \text{id}_{F_{x_{3}}}}\longrightarrow F_{x_{1}} \circ F_{x_{2}} \circ F_{x_{3}} 
\]
and 

\[F_{x_{1}+x_{2}+x_{3}} \stackrel{N_{x_{1},x_{2}+x_{3}}}\longrightarrow F_{x_{1}+x_{2}} \circ F_{x_{3}} 
\stackrel{\text{id}_{F_{x_{1}}} \circ N_{x_{2},x_{3}}}\longrightarrow F_{x_{1}} \circ F_{x_{2}} \circ F_{x_{3}}. 
\]
 The discrepancy between these two compositions defines a constant and there is an additional coherence relation involving $4$ elements of $D$.  This defines a class in $H^{3}(D,\mathbb{C}^{\times})$ where the $\mathbb{C}^{\times}$ appears as the group of automorphisms of the identity functor (or any autoequivalence).  
 
  The last and strongest notion of group action is when the invertible natural transformations $N$ can be chosen so that Diagram (\ref{important}) commutes on the nose.  This can be done if the second obstruction vanishes.

\begin{thm}\label{exist}
 The three different types of representations described above exist for various discrete groups $D$ as described in Assumption \ref{ass}.  Whenever the first obstruction vanishes, the pullback of the above element in $H^{3}(D,\mathbb{C}^{\times})$ to $H^{3}(\pi^{-1}D,\mathbb{C}^{\times})$  is the cohomology class represented by
\[(w_{1},w_{2},w_{3}) \mapsto \Theta^{w_{2}, w_{3}}(w_{1}).\]
where $\Theta^{w_{2}, w_{3}}(w_{1})$ was written explicitly in Equation (\ref{EvalTheta}).\end{thm}
{\bf Proof.}
The stack of $\mathcal{O}_{X}$-modules on $X$ has a canonical equivariant structure for every discrete group action on $X$, it consists for every group element of a map of stacks from the pullback of the stack to itself, and natural transformations for every pair of group elements which satisfy the coherence condition for three group elements.
We consider sheaves on $\mathfrak{G}$ as certain maps of stacks from $\mathfrak{G}$ to the stack of $\mathcal{O}_{X}$-modules on $X$.  Combining this with the above structure, any automorphism of $X$ defines a pullback map from sheaves on $\mathfrak{G}$ to sheaves on the pullback of $\mathfrak{G}$.  If $\mathfrak{G}$ is invariant under the group, we can then chose isomorphisms from the pullback of $\mathfrak{G}$ by any group element to $\mathfrak{G}$.  This gives us for any group element an autoequivalence of abelian categories of $\text{Mod}(\mathcal{O}_{\mathfrak{G}}).$
The functors $F_{x}$ are defined by 
\[F_{x}\mathcal{S} = \mathcal{M}^{x} \otimes t_{x}^{*} \mathcal{S}
\]
where $\mathcal{M}^{x}$ is a trivialization of $\mathfrak{G} \otimes (x \cdot \mathfrak{G})^{-1}$.
These $\mathcal{M}^{x}$ are the twisted line bundles corresponding to $\tau^{w}$ for some $w$ lifting $x$ and were defined in Equation (\ref{eqn:tauDef}).   The cocycle $\Xi^{w_{1}, w_{2}}$ defined in Equation (\ref{xidef}) describes the line bundles $\mathcal{L}^{x_{1},x_{2}}$.
The natural transformations $N_{x_{1}, x_{2}}$ correspond to the trivializations $\theta^{w_{1}, w_{2}}$ defined in Equation (\ref{ThetaEqn}) for some lifts $w_{1},w_{2}$ of $x_{1}, x_{2}$.
\ \hfill $\Box$
\begin{cor}
We make the same assumptions here as in Theorem \ref{exist}.  Suppose also that the first obstruction vanishes.  Then the gerbal representation from the above theorem gives rise to a class in $H^{3}(\pi^{-1}D,\mathbb{C}^{\times})$.   In the case that $\pi^{-1}D \subset \Lambda(\mathfrak{G})_{\mathbb{Z}},$ this class is represented by the cocycle
\[(w_{1}, w_{2}, w_{3}) \mapsto \exp(-\frac{3}{2}E(w_{1}, w_{2}, w_{3})).\]  In the  case that 
$\pi^{-1}D \subset \Lambda(\mathfrak{G})_{(1,1)},$ this class is represented by the cocycle  
\[(w_{1},w_{2},w_{3}) \mapsto \exp(6E(w_{1}, w_{2}, w_{3})).\] 
\end{cor}
{\bf Proof.}
This follows from the discussion in Section \ref{SecondObs}.
\ \hfill $\Box$

When this second obstruction is non-trivial one knows \cite{FZ} that for no choice of functors and natural transformations as above will diagram (\ref{important}) commute.  

\begin{rmk}The above three types of actions also have local versions for which the group acts on the $\mathcal{O}$-linear stack $\mathcal{M}od(\mathcal{O}_{\mathfrak{G}})$ of 
$\mathcal{O}_{\mathfrak{G}}-$modules on $X$.  Here, the functors $F_{x}$ are replaced by maps of stacks from the pullback of the stack to itself.  We leave it to the interested reader to write that out.  This is important in that we would like also to define representations on ``cohomology categories" which are defined for an $\mathcal{O}$-linear stack, whereas here we have only discussed the representation on the zeroth cohomology given by the global sections.
\end{rmk}

%%%%%%%%%%%%%%ADDRESSES%%%%%%%%%%%%%%%%%%
\vskip 0.2in \noindent {\scriptsize {\bf Oren Ben-Bassat,}
Department of Mathematics, University of Haifa, Haifa, Israel, ben-bassat@math.haifa.ac.il}

%%%%%%%%%%%%%%%%%%%%%%%%%%%%%%%%%%%%%%%%%%

\end{document}